\documentclass[11pt,leqno]{amsart}
\usepackage{amsmath,amsfonts,amssymb,amscd,amsthm,amsbsy,upref}
\numberwithin{equation}{section}
\newtheorem{thm}{Theorem}[section]
\newtheorem{lem}[thm]{Lemma}

\newtheorem*{THM}{Theorem}
\newtheorem{cor}[thm]{Corollary}
\newtheorem{prop}[thm]{Proposition}
\theoremstyle{definition}
\newtheorem{defin}[thm]{Definition}

\newtheorem{remark}[thm]{Remark}
\newtheorem{remarks}[thm]{Remarks}
\newtheorem{prob}[thm]{Problem}
\def\nat{{\mathbb N}}
\def\real{{\mathbb R}}
\def\U{{\mathcal U}}
\def\ep{\varepsilon}
\def\supp{\operatorname{supp}}
\def\llbrack{\lbrack\!\lbrack}
\def\rrbrack{\rbrack\!\rbrack}
\begin{document}
\title{On Asymptotically Symmetric Banach Spaces}
\author{M. Junge, D. Kutzarova  \and E. Odell}
\address{M. Junge\\ Department of Mathematics\\ University of Illinois at
Urbana-Champaign,  Urbana, IL 61801}
\email{mjunge@math.uiuc.edu}
\address{D. Kutzarova\\Institute of Mathematics, Bulgarian Academy of
Sciences, Sofia, Bulgaria\\
Current address:
Department of Mathematics\\ University of Illinois at
Urbana-Champaign, Urbana, IL 61801}
\email{denka@math.uiuc.edu}
\address{Department of Mathematics\\ The University of Texas at Austin\\
1 University Station C1200\\ Austin, TX 78712-0257}
\email{odell@math.utexas.edu}

\thanks{Research of the first and third
named authors was partially supported by the National Science Foundation.}

\begin{abstract}
A Banach space $X$ is {\em asymptotically symmetric}  (a.s.)
if for some $C<\infty$, for all $m\in\nat$, for all bounded sequences
$(x_j^i)_{j=1}^\infty \subseteq X$, $1\le i\le m$, for all permutations
$\sigma$ of $\{1,\ldots,m\}$ and all ultrafilters $\U_1,\ldots,\U_m$ on $\nat$,
$$\lim_{n_1,\U_1} \ldots \lim_{n_m,\U_m} \bigg\| \sum_{i=1}^m x_{n_i}^i\bigg\|
\le C\lim_{n_{\sigma (1)},\U_{\sigma (1)}} \ldots
\lim_{n_{\sigma(m)},\U_{\sigma(m)}}
\bigg\|\sum_{i=1}^m x_{n_i}^i\bigg\| .$$
We investigate a.s.\ Banach spaces and several natural variations.
$X$ is weakly a.s. (w.a.s.) if the defining condition holds when
restricted to weakly convergent sequences $(x_j^i)_{j=1}^\infty$.
$X$ is w.n.a.s.\ if we restrict the condition further to normalized
weakly null sequences.

If $X$ is a.s.\ then all spreading models of $X$ are uniformly symmetric.
We show that the converse fails.
We also show that w.a.s.\ and w.n.a.s.\ are not equivalent properties and that
Schlumprecht's space $S$ fails to be w.n.a.s.
We show that if $X$ is separable and has the property that every normalized
weakly null sequence in $X$ has a subsequence equivalent to the unit vector
basis of $c_0$ then $X$ is w.a.s..
We obtain an analogous result if $c_0$ is replaced by $\ell_1$ and also
show it is false if $c_0$ is replaced by $\ell_p$, $1<p<\infty$.

We prove that if $1\le p<\infty$ and $\|\sum_{i=1}^n x_i\|\sim n^{1/p}$
for all $(x_i)_1^n\in \{X\}_n$, the $n^{th}$ asymptotic structure of $X$,
then $X$ contains an asymptotic $\ell_p$, hence w.a.s.\ subspace.
\end{abstract}

\dedicatory{Dedicated to Haskell Rosenthal on the occasion of his 65th
birthday.}

\maketitle

\setcounter{section}{-1}
\section{Introduction}

In their fundamental paper
\cite{KM} J.L.~Krivine and B.~Maurey
introduced
the notion of a stable Banach space and proved that such spaces must contain
almost isometric copies of $\ell_p$ for some $1\le p<\infty$.
$X$ is {\em stable\/} means that for all ultrafilters $\U_1$ and $\U_2$ on
$\nat$ and all bounded sequences $(x_n)$ and $(y_n)$ in $X$,
\begin{equation}\label{eq:intro1}
\lim_{n,\U_1} \lim_{m,\U_2} \|x_n + y_m\| = \lim_{m,\U_2} \lim_{n,\U_1}
\|x_n + y_m\|\ .
\end{equation}

A major application of \cite{KM} was to deduce Aldous' theorem \cite{A}
that every infinite dimensional subspace of $L_p$ ($1\le p<2$) contains
almost isometric copies of $\ell_r$ for some $p\le r\le 2$, by proving
that $L_p$ is stable.
The starting point of our investigation is a result on noncommutative
$L_p$ spaces.
We first formulate the problem which led to that research.

\begin{prob}\label{prob1}
Let $1\le p\le2$ and let $X$ be an infinite dimensional subspace of
$L_p(N)$, the noncommutative $L_p$ space associated with a von~Neumann
algebra $N$.
Does $X$ contain an isomorph of $\ell_r$ for some $p\le r\le 2$?
\end{prob}

In general \cite{JR} noncommutative $L_p$ spaces fail to be stable.
However, if $1<p<\infty$ these spaces $L_p(N)$ satisfy the following
inequality for some universal $C<\infty$ \cite{JR}.
For all $m$, all permutations $\sigma$ of $\{1,\ldots,m\}$, all bounded
sequences in $L_p(N)$, $(x_n^i)_{n=1}^\infty$ for $i\le m$, and all
ultrafilters $\U_1,\ldots,\U_m$ on $\nat$,
\begin{equation}\label{eq:intro2}
\lim_{n_1,\U_1}\ldots \lim_{n_m,\U_m} \bigg\| \sum_{i=1}^m x_{n_i}^i\bigg\|
\le C\lim_{n_{\sigma (1)},\U_{\sigma (1)}} \ldots
\lim_{n_{\sigma(m)},\U_{\sigma (m)}}
\bigg\| \sum_{i=1}^m x_{n_i}^i\bigg\| \ .
\end{equation}
Every stable space $X$ satisfies \eqref{eq:intro2} with $C=1$ and if $Y$
is isomorphic to a stable space, then $Y$ satisfies \eqref{eq:intro2}
for some $C$.

$X$  is called {\em asymptotically symmetric\/} (a.s.) if it satisfies
\eqref{eq:intro2} for some $C<\infty$ where the sequences $(x_n^i)\subseteq X$
\cite{JR}.
It is easy to check that the Tsirelson space $T$ is a.s.\ and so a.s.\ need
not imply that a space contains an isomorph of some $\ell_p$ or $c_0$.
But the following problem remains open.

\begin{prob}\label{prob2}
Let $X$ be an infinite dimensional asymptotically symmetric Banach space.
Does $X$ contain an asymptotic $\ell_p$ subspace for some $p$?
\end{prob}

A Banach space $X$ with a basis $(e_i)$ is called {\em asymptotically $\ell_p$}
 \cite{MT} if for some $C<\infty$ for all $n$, every normalized block basis
$(x_i)_{i=1}^n$ of $(e_i)_{i=n}^\infty$ is $C$-equivalent to the
unit vector basis of $\ell_p^n$. In this case $(e_i)$ is called an
{\em asymptotic $\ell_p$ basis} for $X$.

If $X$ is a.s.\ then every spreading model $(\tilde x_i)$ of a normalized
basic sequence $(x_i)$ in $X$ is $C$-symmetric for some fixed $C<\infty$.
As we shall see in \S2, a.s.\ says much more than this.

\begin{thm}\label{thm:intro3}
There is a reflexive Banach space $X$ such that all spreading models of $X$
are $C$-symmetric for some fixed $C$ yet $X$ is not asymptotically
symmetric.
\end{thm}

In fact given $1<p<\infty$ we can construct the space $X$ in
Theorem~\ref{thm:intro3} to have the property that every normalized
weakly null sequence in $X$ has a subsequence 4-equivalent to the unit
vector basis of $\ell_p$.
We show that if $\ell_p$ is replaced by $c_0$ (or $\ell_1$
under an appropriate restatement) then one has positive results.

In \S3 we prove

\begin{thm}\label{thm:intro4}
There exists a reflexive space $Y$ which is not asymptotically symmetric
and yet for some $C<\infty$, \eqref{eq:intro2} holds for all
normalized weakly null sequences $(x_n^i)_{n=1}^\infty$, $i\le m$, in $Y$.
\end{thm}

$Y$ is the Tirliman space $T_i$ $(2;1/2)$, a subsymmetric version of a
space due to L.~Tzafriri \cite{Tz}, as presented in \cite{CS}.
We show that $Y$ contains an asymptotic $\ell_2$ subspace and hence
is not minimal.
More generally, we prove that if $X$ satisfies, $K^{-1}n^{1/p} \le \|
\sum_{i=1}^n x_i\| \le Kn^{1/p}$ for some $1\le p<\infty$, $K<\infty$,
all $n\in\nat$ and all $(x_i)_{i=1}^n \in \{X\}_n$, then $X$ contains an
asymptotic $\ell_p$, hence w.a.s., subspace.
$\{X\}_n$ is the $n^{th}$ asymptotic structure of a space $X$
due to \cite{MMT}.

In \S4 we show that Tsirelson's space $T$ is not iteration stable.
This notion, due to H.~Rosenthal, is another weakening of the definition of
stability.

\S1 contains the definitions of certain variants of a.s. and the relations
between them as well as certain preliminaries.
The authors wish to thank H.~Rosenthal for many enlightening discussions.

\section{Preliminaries}

The definition of an asymptotically symmetric Banach space $X$ can also be
formulated in this way.
$X$ is {\em asymptotically symmetric\/} if for some $C<\infty$, for all
$m\in \nat$, for all bounded sequences $(x_n^i)_{n=1}^\infty \subseteq X$,
$1\le i\le m$, and for all permutations $\sigma$ of $\{1,2,\ldots,m\}$,
\begin{equation}\label{prelim1}
\lim_{n_1\to\infty} \ldots \lim_{n_m\to\infty}
\bigg\|\sum_{i=1}^m x_{n_i}^i \bigg\|
\le C\lim_{n_{\sigma(1)}\to\infty} \ldots \lim_{n_{\sigma (m)}\to\infty}
\bigg\| \sum_{i=1}^m x_{n_i}^i\bigg\|\ ,
\end{equation}
provided that these iterated limits exist.

Just as stability was weakened to ``weak stability'' to handle the case
of spaces like $c_0$ \cite{ANZ}, it is useful to consider certain variants
of a.s.
We will say that $X$ is {\em weakly asymptotically symmetric\/} (w.a.s.) if
\eqref{prelim1} holds when restricted to sequences $(x_n^i)_{n=1}^\infty
\subseteq X$ which are weakly null.
This is actually equivalent to restricting to weakly convergent sequences,
although the constant $C$ could vary.
One could also restrict \eqref{prelim1} to normalized sequences in $X$.
However it is not hard to show that this would be the same as a.s.
$X$ is said to be {\em weakly null normalized asymptotically symmetric\/}
(w.n.a.s.) if \eqref{prelim1} holds for all normalized weakly null
sequences $(x_n^i)_{n=1}^\infty \subseteq X$.
In \S3 we will show that w.a.s. and w.n.a.s. are not equivalent properties.

If $X$ has a basis $(e_i)$ then $X$ is {\em block asymptotically symmetric\/}
(b.a.s.) with respect to $(e_i)$ if \eqref{prelim1} holds for all
bounded block bases $(x_n^i)_{n=1}^\infty$, $i\le m$, of $(e_i)$.

Recall that $(\tilde x_i)$ is a spreading model of a normalized basic
sequence $(x_i)$ if for all $\ep >0$ and $k\in \nat$ there exists $n\in\nat$
so that for all $n\le i_1<\cdots < i_k$ and scalars $(a_i)_1^k$ with each
$|a_i| \le 1$,
\begin{equation}\label{prelim2}
\left|\bigg\| \sum_{i=1}^k a_i x_{n_i}\bigg\| -
\bigg\| \sum_{i=1}^k a_i \tilde x_i\bigg\| \right| <\ep\ .
\end{equation}
Using  Ramsey's theorem it is easy to show that every normalized basic
sequence admits a subsequence $(x_i)$ generating a spreading model (for more
on spreading models see \cite{BL}).

If $(\tilde x_i)$ is a spreading model of $(x_i)$ then $(\tilde x_i)$ is
1-subsymmetric, i.e.,
$$\Big\| \sum_{i=1}^n a_i \tilde x_i\Big\| =
\Big\| \sum_{i=1}^n a_i \tilde x_{k_i}\Big\|$$
for all $(a_i)_1^n$ and $k_1<\cdots < k_n$.
(Some authors call this 1-spreading and reserve the subsymmetric definition
for unconditional spreading sequences.)
The fact that if $X$ is a.s. with constant $C$ implies every spreading model
$(\tilde x_i)$ of $X$ generated by $(x_i)$ is {\em $C$-symmetric\/}, i.e.,
$\|\sum_{i=1}^m a_i\tilde x_i\|\le C\|\sum_{i=1}^m a_i\tilde x_{\sigma (i)}\|$
for all permutations $\sigma$ of $\nat$, follows easily from \eqref{prelim1}
setting $x_n^i = a_i x_n$ for all $i,n$.
(Some authors call this $C$-exchangeable and define the $C$-symmetric constant
via $\|\sum_{i=1}^m \pm a_i \tilde x_i\|$, over all choices of signs.
A $C$-symmetric basis is unconditional, but not necessarily $C$-unconditional.)

There are several known means of joining the infinite and finite dimensional
structure of a Banach space in a conceptual manner.
In addition to spreading models, we have the theory of asymptotic
structure \cite{MMT} (see \S3), which, in particular, gives rise to
the definition of asymptotic $\ell_p$ bases presented above.
A similar definition can be made for {\em asymptotically $c_0$ bases\/}.
A third concept is an {\em asymptotic model\/} of $X$ \cite{HO}.
These are generated much like spreading models except one uses a certain
infinite array of normalized basic sequences in $X$,
$(x_j^i)_{j=1}^\infty$, $i\in\nat$, and one replaces
$\| \sum_{i=1}^k a_i x_{n_i}\|$ in \eqref{prelim2} by
$\| \sum_{i=1}^k a_i x_{n_i}^i\|$.
The asymptotically symmetric definition and variants constitute a fourth
way to join the finite and infinite geometries of a space.

We gather together some of the relationships between these four concepts
in our next proposition.
For the sake of completeness we include some of our above observations.
The ``moreover'' statement requires two of our main results below.

\begin{prop}\label{prop:prelim1}
Let $X$ be an infinite dimensional Banach space.
\begin{itemize}
\item[a)] If $X$ is a.s. then the spreading models of $X$  are uniformly
symmetric.
Similarly, if $X$ is b.a.s. with respect to $(e_i)$ then the spreading models
of  block  bases of $(e_i)$ are uniformly symmetric.
\item[b)] $X$ is a.s. $\Rightarrow$ ($X$ is b.a.s. if $X$ has a basis)
$\Rightarrow X$ is w.a.s. $\Rightarrow X$ is w.n.a.s.
\item[c)] If $X$ is reflexive, then $X$ is a.s. $\Leftrightarrow X$ is
w.a.s $\Leftrightarrow$ ($X$ is b.a.s., if $X$  has a basis).
\item[d)] If $(e_i)$ is a boundedly complete basis for $X$ then $X$ is
a.s. $\Leftrightarrow X$ is b.a.s. with respect to $(e_i)$.
\item[e)] If $X$ is a.s. and not reflexive, then $\ell_1$ is isomorphic to a
spreading model of $X$.
\item[f)] If all asymptotic models of $X$ are symmetric then $X$ is w.a.s.
(and b.a.s if it has a basis).
\item[g)] If $X$ has an asymptotically $\ell_p$ or $c_0$ basis $(e_i)$,
then $X$ is b.a.s. with respect to $(e_i)$.
\end{itemize}

Moreover all of the converses of the one sided implications in
a)--g) are false, in general.
\end{prop}

\begin{proof}
a)--c) follows easily from our previous remarks and d) which holds by
standard gliding hump arguments.
To see e)  we note that a nonreflexive $X$ contains a non weakly null
normalized basic sequence $(x_i)$ having a spreading model $(\tilde x_i)$
which thus satisfies for some $\delta >0$, $\|\sum a_i\tilde x_i\|\ge
\delta \sum a_i$ if the $a_i$'s are nonnegative.
If $X$ is a.s. then $(\tilde x_i)$ is symmetric, hence unconditional,
hence equivalent to the unit vector basis of $\ell_1$.

f) If all asymptotic models of $X$ are symmetric, then it follows easily
from Krivine's theorem \cite{K} that for some $C<\infty$ if $(x_i)$ is an
asymptotic model of $X$ generated by either a weakly null array or a block
basis array (if $X$ has a basis), then $(x_i)$ is $C$-equivalent to the
unit vector basis of $c_0$ or $\ell_p$ for some fixed (independent of $(x_i)$)
$p\in [1,\infty)$ \cite{HO}.
Thus $X$ is w.a.s. (or b.a.s.) by much the same easy argument that
yields g).

The converses in a) are shown to fail in \S2.
By considering the summing basis for $c_0$, which is subsymmetric but not
symmetric, we see that the second converse in b) fails.
$c_0$ also provides a counterexample to the first converse since it is not
a.s.\ but it is b.a.s.\ with respect to the unit vector basis.
The third converse is proved false in \S3.
Every normalized unconditional basic sequence in $L_p$ is equivalent to
an asymptotic model of $L_p$, $1<p<\infty$ \cite{HO} and thus the
converses in f) also fail.
$L_p$ also provides a converse to g).
\end{proof}

One can also ask variations of problem~\ref{prob2} using these alternate
asymptotic notions.

\begin{prob}\label{prob:prelim2}
Assume that for some $C<\infty$ and $1\le p<\infty$, all spreading models
of $X$ (or even all asymptotic models of $X$) are $C$-equivalent to the
unit vector basis of $\ell_p$.
Does $X$ contain an asymptotic $\ell_p$ basic sequence?
\end{prob}

\section{Spreading models and asymptotic symmetry}            %% section 2

The following theorem easily yields Theorem~\ref{thm:intro3}.

\begin{thm}\label{thm:2.1}
Let $1<p<\infty$.
There exists a reflexive infinite dimensional Banach space $X$ which is
not asymptotically symmetric and yet satisfies the following.
\begin{equation}\label{eq:2.1}
\begin{cases}
\text{Every normalized basic sequence in $X$ admits a subsequence}\cr
\text{that is 4-equivalent to the unit vector basis of $\ell_p$}
\end{cases}
\end{equation}
\end{thm}

\begin{proof}
We shall define spaces $X_k$ for $k\in \nat$ and set
$X=(\sum_{k=1}^\infty X_k)_p$.
Each $X_k$ will be reflexive and satisfy \eqref{eq:2.1} with 4 replaced
by $3+\ep$ for any $\ep>0$.
Gliding hump arguments then yield \eqref{eq:2.1} for $X$.

Fix $1<q<p<r$ with $\frac1q +\frac1r =1$ and let $k\in\nat$.
$X_k$ will have a normalized 1-unconditional basis $\{e_j^i :1\le i\le k$,
$j\ge i\}$ which we visualize as $k$ infinite rows of an upper triangle array.
We will define the norm on $X_k$ so that if $n_1<\cdots < n_k$ then
$\|\sum_{i=1}^k e_{n_i}^i\| = k^{1/r}$  and we shall say that the collection
$(e_{n_i}^i)_{i=1}^k$ is {\em permissible\/}.
In addition we will have that $\|e_{n_k}^1 +\cdots + e_{n_1}^k\|=1$.
Thus
$$\lim_{n_1\to\infty}\ldots\lim_{n_k\to\infty} \|e_{n_1}^1+\cdots +e_{n_k}^k\|
= k^{1/r} \lim_{n_k\to\infty} \ldots \lim_{n_1\to\infty}
\|e_{n_1}^1 +\cdots + e_{n_k}^k\|\ .$$
Hence $X$ is not a.s.

Let $(a_i)_{i=1}^k\in B_{\ell_q^k}$, the unit ball of $\ell_q^k$.
Let $x\in \text{span}\{e_j^i : i\le k$, $j\ge i\}$, say
$x= \sum b_j^i e_j^i$.
We define
$$|x|_{(a_i)_{i=1}^k} = \sup
\left( \sum_{j=1}^\infty \bigg(\sum_{i=1}^k |a_ib_{n_j^i}^i|\bigg)^p
\right)^{1/p}$$
where the ``sup'' is taken over all lexicographically ordered integers
$n_1^1 <n_1^2 <\cdots < n_1^k < n_2^1 < n_2^2 <\cdots$.
Thus each collection $(b_{n_j^i}^i)_{i=1}^k$ are the coordinates of $x$ with
respect to the permissible collection $(e_{n_j^i}^i)_{i=1}^k$ and these
collections move to the right in our array picture as $j$ increases.

$X_k$ is the completion of
$(\text{span}\{e_j^i : i\le k$, $j\ge i\},\|\cdot\|)$
where
$$\|x\| = \sup \{ |x|_{(a_i)_{i=1}^k} : (a_i)_{i=1}^k \in B_{\ell_q^k}\}\ .$$
If $(e_{n_i}^i)_{i=1}^k$ is a permissible collection and
$x= \sum_{i=1}^k b_i e_{n_i}^i$, then
$\|x\| = (\sum_{i=1}^k |b_i|^r)^{1/r}$.
Indeed the lower estimate is immediate and suppose that
$\|x\| =  |x|_{(a_i)_{i=1}^k}$ for some $(a_i)_{i=1}^k\in B_{\ell_q^k}$.
We may thus write
$$\|x\| = \left[ \bigg( \sum_{i\in I_1} |a_i b_i|\bigg)^p
+\cdots +
\bigg( \sum_{i\in I_t} |a_i b_i|\bigg)^p
\right]^{1/p}$$
where $I_1,\ldots, I_t$ are disjoint subsets of $\{1,\ldots, k\}$.
Thus, since $\frac1q + \frac1r =1$,
$$\|x\| \le \sum_{i=1}^k |a_i b_i| \le \|(b_i)_{i=1}^k \|_r\ .$$

Next we let $y = \sum_{i=1}^k e_{n_i}^i$ where  $n_k <n_{k-1} < \cdots <n_1$.
Let $\|y\| = |y|_{(a_i)_{i=1}^k}$.
Since any permissible collection of $e_j^i$'s will intersect the support
of $y$ in at most one coordinate we have that
$$|y|_{(a_i)_{i=1}^k}
\le \bigg( \sum_{i=1}^k |a_i|^p\bigg)^{1/p}
\le \| (a_i)_1^k \|_q \le 1\ .$$
This completes the proof of our assertions which yield that $X$ is not a.s.

The  basis for $X_k$ is boundedly complete and so once we prove that $X_k$
satisfies \eqref{eq:2.1} with constant $3$
and for all normalized block bases it
will follow that $X_k$ is reflexive and satisfies \eqref{eq:2.1} for
$3+\ep$ and for all normalized basic sequences.

Let $\ep>0$ and let $(x_i)$ be a normalized block basis of $\{e_j^i: i\le k$,
$j\ge i\}$.
Passing to a subsequence we may assume that for all $m$
\begin{equation}\label{eq:2.2}
\max \{j:e_j^i \in\supp x_m\} + 2k
< \min \{j:e_j^i \in \supp  x_{m+1}\}
\end{equation}
and for some $(a_i)_{i=1}^k \in B_{\ell_q^k}$,
\begin{equation}\label{eq:2.3}
|x_m|_{(a_i)_{i=1}^k} > 1-\ep\ .
\end{equation}
Since \eqref{eq:2.2} spaces the $x_m$'s to have at least $2k$ ``columns''
between successive supports we obtain, using \eqref{eq:2.3}, that for
all scalars $(c_m)$,
$$\bigg\|\sum_m c_m x_m\bigg\|\ge(1-\ep)\bigg(\sum_m |c_m|^p\bigg)^{1/p}\ .$$
Indeed one can string together the lexicographically ordered lists that
yield each norm $|x_m|_{(a_i)_1^k}$, inserting extra elements as
needed into the gaps.

It remains to prove that
$\|\sum_m c_m x_m\| \le 3(\sum_m |c_m|^p)^{1/p}$.
Let $x= \sum_m c_m x_m$ and suppose $\|x\| = |x|_{(a_i)_1^k}$.
We let $x= \sum b_j^i e_j^i$ and choose $n_1^1 <\cdots < n_1^k <n_2^1<n_2^2
< \cdots$ so that
$$\|x\| = \bigg( \sum_{j=1}^\infty \bigg( \sum_{i=1}^k |a_ib_{n_j^i}^i|
\bigg)^p \bigg)^{1/p}\ .$$
For each $j$ let $A_j = \{e_{n_j^i}^i\}_{i=1}^k$ be the corresponding
permissible collection.
\begin{itemize}
\item[{}] Let $J_1 = \{j:A_j$ intersects the support of exactly one $x_m\}$.
\item[{}] Let $J_0 = \{j:A_j$ intersects the support of more than one $x_m\}$.
\end{itemize}
and let $J_0 = J_2 \cup J_3$ where $J_2$ contains every other integer in $J_0$
and $J_3 = J_0\setminus J_2$.
Thus if $j_1$ and $j_2$ are distinct integers in $J_2$ (or $J_3$) then
$A_{j_1}$ and $A_{j_2}$ cannot both intersect the support of the same $x_m$.

By the triangle inequality in $\ell_p$,
$$\|x\| \le \sum_{\ell=1}^3 \bigg( \sum_{j\in J_\ell}
\bigg( \sum_{i=1}^k |a_i b_{n_j^i}^i|\bigg)^p \bigg)^{1/p}\ .$$
We shall show that each of these three terms is bounded above by
$\| (c_m)\|_p$.

For $m \in \nat$ let $I_m = \{j\in J_1 :A_j\cap \supp x_m \ne\emptyset\}$.
Thus $I_m \cap I_{m'} = \emptyset$ if $m \ne m'$.
Since $|c_m| = \|c_m x_m\|$ we have
$$\bigg( \sum_{j\in J_1} \bigg( \sum_{i=1}^k |a_i b_{n_j^i}^i|\bigg)^p
\bigg)^{1/p} =
\bigg( \sum_m \sum_{j\in I_m} \bigg( \sum_{i=1}^k |a_i b_{n_j^i}^i|\bigg)^p
\bigg)^{1/p}
\le \|(c_m)\|_p .$$

Finally we estimate the $J_2$ sum (the $J_3$ estimate is identical).
The $J_2$ sum is an $\ell_p$ sum of terms of the form
$$Q_j = \sum_{t=1}^s |c_{m_t}| \sum_{i\in I_t} |a_i d_{t,i}|$$
where $I_1 <\cdots < I_s$ are subsets of $\{1,\ldots,k\}$ and
$m_1<\cdots < m_s$.
$x_{m_1},\ldots ,x_{m_s}$ are those $x_i$'s for which $A_j\cap \supp
x_i\ne \emptyset$ and the $d_{t,i}$'s are the corresponding relevant
coordinates of $x_{m_t}$.
$(x_{m_i})_{i=1}^s$ depends of course upon $j$ but these families are
disjoint for different $j$'s $\in J_2$.
Thus it suffices to prove that
$$Q_j \le \bigg( \sum_{t=1}^s |c_{m_t}|^p\bigg)^{1/p}\ .$$
Now for $t$ fixed $(d_{t,i})_{i\in I_t}$ are the coordinates of $x_{m_t}$
with respect to a subset (indexed by $I_t$) of a permissible collection of
$e_v^u$'s
and these are in turn 1-equivalent to the unit vector basis of
$\ell_r^{|I_t|}$ as we have shown.
Hence
\begin{equation*}
\begin{split}
Q_j &\le \sum_{t=1}^s |c_{m_t}| \, \|(a_i)_{i\in I_t}\|_q\\
&\le \|(c_{m_t})_{t=1}^s\|_r \|(a_i)_{i=1}^k\|_q\\
&\le \|(c_{m_t})\|_r \le \|(c_{m_t})\|\ .
\end{split}
\end{equation*}
\end{proof}

\begin{remarks}\label{rem:2.2}
The space $X$, constructed in Theorem~\ref{thm:2.1}, satisfies \eqref{eq:2.1}
and also has the property that $\ell_r$ and $c_0$ are asymptotic versions of
$X$ (see \cite{MMT}).
A different example of this sort of phenomenon is given in \cite{OS} where
a reflexive space $Z$ is constructed satisfying \eqref{eq:2.1} with 4 replaced
by $1+\ep$, $\ep>0$ arbitrary, yet $Z$ has $\ell_r$ as an asymptotic
version for some $r\ne p$.
It would be interesting to see if one could construct $X$ as in
Theorem~\ref{thm:2.1} to satisfy \eqref{eq:2.1} for 4 replaced by $1+\ep$.
Another natural question is to ascertain what happens if $\ell_p$ is
replaced by $\ell_1$ or $c_0$.
We will show that in these cases one obtains positive results.

%$X$ is called $K$-strong Schur if for all normalized $(x_i)\subseteq X$ with
%$\|x_i - x_j\| \ge \delta >0$ for $i\ne j$, some subsequence of $(x_i)$
%is $K/\delta$-equivalent to the unit vector basis of $\ell_1$.
%For example $\ell_1$ is 2-strong Schur under this quantification.
\end{remarks}

\begin{thm}\label{thm:2.3}
%Let $X$ be a $K$-strong Schur space with basis $(e_i)$.
Let $X$ have a basis $(e_i)$.
Assume that for some $K<\infty$ every spreading model of any normalized
block basis of $(e_i)$ is $K$-equivalent to the unit vector basis of $\ell_1$.
Then $X$ is block asymptotically symmetric with respect to $(e_i)$.
\end{thm}

\begin{proof}
By renorming $X$ we may assume that $(e_i)$ is a bimonotone basis.
Let $m\in\nat$, $(b_i)_1^m \subseteq [-1,1]$,
$(x_j^i)_{j=1}^\infty$ a normalized block basis of $(e_i)$ for $i\le m$ and
$\sigma$ a permutation of $\{1,\ldots,m\}$ so that the iterated limits
$$\lim_{n_1\to\infty}\ldots \lim_{n_m\to\infty} \bigg\| \sum_{i=1}^m
b_i x_{n_i}^i\bigg\| =1\ \text{ and }\
\lim_{n_{\sigma (1)}\to\infty} \ldots \lim_{n_{\sigma (m)}\to\infty}
\bigg\| \sum_{i=1}^m b_i x_{n_i}^i\bigg\|$$
both exist, the first being equal to 1.

We visualize $(x_j^i)_{i,j}$ as an array of $m$ infinite rows.
Using  Ramsey's theorem, by passing to a subsequence of the columns, given
$\ep>0$ we may assume that for all integers $n_1<\cdots < n_m$ and
$k_1< \cdots <k_m$, for all $f\in B_{X^*}$ there exists $g\in B_{X^*}$ with
$|f(x_{n_i}^i) - g(x_{k_i}^i)| \le \ep$ for $i\le m$.
This follows by partitioning $[-1,1]$ into finitely many intervals
$(I_t)_{t=1}^\ell$ of length less than $\ep$ and thus inducing a finite
coloring of $[\nat]^m$ as follows:
$(n_1,\ldots,n_m)$ has color $(I_{t_1},\ldots,I_{t_m})$ if there exists
$f\in B_{X^*}$ with $f(x_{n_i}^i)\in I_{t_i}$ for $i\le m$.

It follows that up to an arbitrarily small error
$$\Big\| \sum_{i=1}^m b_i x_{n_i}^i \Big\|
\approx \Big\|\sum_{i=1}^m b_i x_{k_i}^i\Big\| \approx 1$$
whenever $n_1<\cdots <n_m$ and $k_1 < \cdots < k_m$.
To avoid trivial but tedious error estimates in the remainder of the proof
we shall assume $\ep =0$.

\noindent
We may also assume that if $x\in \text{span}(x_{j_1}^i)_{i=1}^m$,
$y\in\text{span}(x_{j_2}^i)_{i=1}^m$ for $j_1< j_2$ then
$\supp x< \supp y$ w.r.t. $(e_n)$.
Finally we assume similar stabilizations for the order induced by $\sigma$.

Let $n_1^1 <\cdots < n_1^m < n_2^1 <\cdots < n_2^m <\cdots$ and set
$y_j = \sum_{i=1}^m b_i x_{n_j^i}^i$ for $j\in\nat$.
$(y_j)$ is a normalized block basis of $(e_i)$ and thus, passing to a
subsequence, we may assume it has a spreading model which is
$K$-equivalent to the unit vector basis
of $\ell_1$.
%Moreover there exists $f\in B_{X^*}$ with $f(y_j)\ge 1/K$ for all $j$.

Hence for all $k$ there exists $F\subseteq \nat$ with $|F|=k$ and
$f\in B_{X^*}$ with $f(y_j)\ge 1/K$ for $j\in F$ (up to an arbitrarily
small error).
Using the pigenhole principle and ignoring small errors we obtain $m$ of the
$y_j$'s which we relabel as $y_1,\ldots,y_m$, $f\in B_{X^*}$ and
%Passing to a subsequence again, ignoring arbitrarily small errors, we may
%assume that there exists
$(a_i)_1^m$ so that $f(x_{n_j^i}^i)=a_i$ for all $j$ and $1\le i,j\le m$.

We shall show that
$$1\le K\Big\| \sum_{i=1}^m b_{\sigma (i)} x_{n_{\sigma (i)}}^{\sigma (i)}
\Big\|$$
provided that $n_{\sigma (1)} <\cdots < n_{\sigma (m)}$.
{From} our stabilizations it suffices to produce
$n_{\sigma (1)} <\cdots < n_{\sigma (m)}$ satisfying this.

We choose $n_{\sigma (1)} <\cdots < n_{\sigma (m)}$ so that each
$x_{n_{\sigma (i)}}^{\sigma (i)}$ is in the support of some $y_j$, $j\le m$.
It follows that $f(x_{n_{\sigma (i)}}^{\sigma (i)}) = a_{\sigma (i)}$
and the claim follows.
\end{proof}

The proof yields that the b.a.s. constant of $X$ is bounded by a function
of $K$ and the basis constant of $(e_i)$.

\begin{thm}\label{thm:2.4}
Let $X$ be a separable Banach space such that every
spreading model of a normalized weakly null
sequence in $X$  is equivalent to the unit vector basis
of $c_0$.
Then $X$ is weakly asymptotically symmetric.
\end{thm}

\begin{proof}
It follows from \cite[Proposition 3.2]{AOST} that for some $C<\infty$
every spreading model of a
normalized weakly null sequence in $X$ is
$C$-equivalent to the unit vector basis of $c_0$.

Let $m\in\nat$ and let $(x_j^i)_{j=1}^\infty$, $i\le  m$ be normalized
weakly null sequences in $X$.
Let $K<\infty$, $(b_i)_1^m\subseteq \real$ and assume that for some
permutation $\sigma$
$$\lim_{n_1\to\infty} \ldots \lim_{n_m\to\infty}
\Big\| \sum_{i=1}^m b_i x_{n_i}^i\Big\| =K$$
and
$$\lim_{n_{\sigma (1)}\to\infty} \ldots \lim_{n_{\sigma (m)}\to\infty}
\Big\| \sum_{i=1}^m b_i x_{n_i}^i\Big\| =1\ .$$
We will prove that $K\le C$ which will complete the proof of the theorem.

As in the proof of the previous theorem, by passing to a subsequence of
the columns and ignoring arbitrarily small errors we may assume that
$\|\sum_{i=1}^m b_i x_{n_i}^i\|=K$ if $n_1<\cdots < n_m$.
Moreover, we may assume that if $f\in B_{X^*}$ with $f(x_{n_i}^i) = a_i$
for $i\le m$ and if $k_1<\cdots < k_m$, then there  exists $g\in B_{X^*}$
with $g(x_{k_i}^i)= a_i$ for $i\le m$.

We shall say that $z$ is $((b_i)_1^m,\sigma)$ distributed if $z=\sum_{i=1}^m
b_{\sigma (i)} x_{n_{\sigma (i)}}^{\sigma (i)}$ for some
$n_{\sigma (1)} <\cdots < n_{\sigma (m)}$ and as above, by Ramsey's Theorem,
we may assume that for such a vector, $\|z\|=1$.
In addition we may assume that if $(z_i)_{i=1}^m$ are all
$((b_i)_1^m,\sigma)$ distributed with
$z_j = \sum_{i=1}^m b_{\sigma (i)} x_{n_{\sigma(i)}^j}^i$
and
$n_{\sigma (m)}^j < n_{\sigma (1)}^{j+1}$ for $j<m$,
then $\|\sum_{j=1}^m z_j\|$ does not depend upon the particular choice of
the $n_{\sigma (i)}^j$'s.
Finally, since the rows $(x_j^i)_{j=1}^\infty$ are weakly null, we can assume
that the coordinates supporting such a sequence $(z_j)_{j=1}^m$,
namely $(x_{n_{\sigma (i)}^j}^i)_{i=1,\, j=1}^{m,m}$ are suppression-1
unconditional.
(This argument is used in \cite{HO} and \cite{AOST}).
Roughly, if one has $f\in B_{X^*}$ and one considers
$z =\sum_{i,j=1}^m a_j^i x_{n_{\sigma (i)}^j}^i$
with its coordinates sufficiently spread out then one can ``slide'' the
coordinates $I$ one wishes to kill, preserving the order, so that $f\approx 0$
on these coordinates.
The new vector $w$, distributed exactly the
same as $z$, hence $\|w\| = \| z\|$, satisfies
$$f(w) = \Big| f\Big( \sum_{i,j\notin I} a_j^i  x_{n_{\sigma(i)}^j}^i\Big)\Big|
\le \|w\| = \|z\|\ .$$

Now let $(z_j)_{j=1}^\infty$ be a ``block basis'' of $(x_j^i)$ with each
$z_j$ having $((b_i)_1^m,\sigma)$ distribution.
$(z_j)$ is normalized weakly null so passing to a subsequence we may assume
it has a spreading model which
is $C$-equivalent to the unit vector basis of $c_0$.
%Consider $\|\sum_{i=1}^m z_i\|$.
In particular, by relabeling, we may assume that $\|\sum_{i=1}^m z_i\|\le C$.
By restricting this vector to a suitable set of coordinates we obtain a
vector equal to $\sum_{i=1}^m b_i x_{n_i}^i$ for some $n_1<\cdots < n_m$.
Thus $C\ge \|\sum_{i=1}^m z_i\|\ge \|\sum_{i=1}^m b_i x_{n_i}^i\|  = K$
and the proof is complete.
%Note that the w.a.s. constant can be taken as $C+\ep$ for all $\ep >0$.
\end{proof}

\section{Variants of a.s. and Tsirelson-like spaces}            %%%% 33333

\begin{thm}\label{thm:3.1}
There exists a reflexive Banach space $Y$ which is w.n.a.s. but not w.a.s.
\end{thm}

$Y$ is Tzafriri's space $Ti(2;1/2)$ (\cite[section X.D]{CS}).
We recall the definition.
Let $c_{00}$ be the linear space of finitely supported sequences of reals.
If $x\in c_{00}$ and $E\subseteq \nat$, we set $Ex(i) = x(i)$ if $i\in E$
and $0$ otherwise.
For sets $E,F\subseteq \nat$, $E<F$ denotes $\max E< \min F$.
$Y$ is the completion of $c_{00}$ under the norm given by the following
implicit equation.
\begin{equation}\label{eq:3.1}
\|x\| = \max \left( \|x\|_\infty , \sup \frac1{2\sqrt{n}} \sum_{i=1}^n
\|E_ix\|\right)
\end{equation}
where the ``sup'' is taken over all $n\in \nat$ and $E_1<\cdots < E_n$.

$Y$ is reflexive and the unit vector basis $(e_i)$ is a normalized
1-unconditional 1-subsymmetric basis for $Y$.
We recall two facts from \cite{CS}.

\noindent {\bf Fact 1.}
For all $x\in Y$, $\|x\| \le \|x\|_2$.

\noindent {\bf Fact 2.}
\cite[lemma X.d.4, p.109]{CS}
If $(u_i)_{i=1}^n$ is a finite block basis of $(e_i)$ then
$$\Big\| \sum_{i=1}^n u_i\Big\| \le \sqrt3 \bigg( \sum_{i=1}^n \|u_i\|^2
\bigg)^{1/2}$$

{From} Fact 2 and \eqref{eq:3.1} we have that if $(y_i)_{i=1}^n$ is a
normalized block basis of $(e_i)$ then
$$\sqrt{n}/2 \le \Big\| \sum_{i=1}^n y_i\Big\| \le 3\sqrt{n}\ .$$
This yields that $Y$ is w.n.a.s. with constant $6$.
That $Y$ is not w.a.s. follows from either the next theorem or
\cite{Sa} (see the remarks below).

\begin{thm}\label{thm:3.2}
$c_0$ is finitely representable in $Y$.
Moreover, for some $C<\infty$ (equivalently, for all $C>1$), for all $n$
there exist disjointly supported (w.r.t. $(e_i)$) normalized vectors
$(x_i)_{i=1}^n$ in $Y$ with $(x_i)_{i=1}^n$ being $C$-equivalent
to the unit vector basis of $\ell_\infty^n$.
\end{thm}

First we shall show how this theorem completes the proof of
Theorem~\ref{thm:3.1}.
$(e_i)_{i=1}^\infty$ is 1-subsymmetric and hence is its own spreading
model.
But Theorem~\ref{thm:3.2} and \eqref{eq:3.1} yield  that $(e_i)$ is not
symmetric, which every spreading model of a w.a.s. reflexive space would be.

Tzafriri \cite{Tz} constructed a symmetric version of $Y$ denoted by
$V_{1/2,2}$.
The definition of the norm is given by \eqref{eq:3.1} where the sup is
taken over disjoint sets $(E_i)_1^n$ in $\nat$.
The space $V_{1/2,2}$ has finite cotype and hence does not contain
$\ell_\infty^n$'s uniformly.
Thus, a consequence of Theorem~\ref{thm:3.2} is the following corollary
which answers a question from \cite{CS}.

\begin{cor}\label{cor:3.3}
The spaces $Y$ and $V_{1/2,2}$ are not isomorphic.
\end{cor}

This question was answered independently by
B.~Sari  \cite{Sa} who used different techniques.
In fact Sari has proved that $Y= Ti(2;1/2)$
does not contain a symmetric basic sequence.
Moreover, Sari's result also yields that $Y$ is not w.a.s. since $(e_i)$
is not symmetric.
We include Theorem~\ref{thm:3.2} because it is of separate interest.

\begin{lem}\label{lem:3.4}
For $n\ge 4$, $\|\sum_{i=1}^n e_i\| = \sqrt{n}/2$.
\end{lem}

\begin{proof}
The lower estimate is immediate and the case $n=4$
is easy using \eqref{eq:3.1}.
For $n>4$ there exist $E_1<\cdots <E_k$ with
$\|\sum_{i=1}^n e_i\| =\frac1{2\sqrt k} \sum_{j=1}^k\|E_j(\sum_{i=1}^n e_i)\|$.
Let $|E_j\cap \{1,\ldots,n\}| = n_j$, hence
$\sum_{j=1}^k n_j \le n$.
By Fact~1 and Cauchy-Schwarz,
\begin{equation*}
\begin{split}
\Big\| \sum_{i=1}^n e_i\Big\|
& \le \frac1{2\sqrt k} \sum_{j=1}^k \sqrt{n_j}
\le \frac1{2\sqrt k} \bigg(\sum_{j=1}^k n_i\bigg)^{1/2}\sqrt{k}\\
&\le \frac{\sqrt{n}}2\ .
\end{split}
\end{equation*}
\end{proof}

\begin{lem}\label{lem:3.5}
Let $\delta = \sqrt3 /2 <1$ and let $(u_i)_1^n$ be a block basis of $(e_i)$
with $\|u_i\|\le 1$ for $i\le n$.
Let $E_1<\cdots < E_k$ be subsets of $\nat$, so that for $i\le n$,
$\supp u_i$ intersects at most one $E_j$.
Then
$$\frac1{2\sqrt k} \sum_{j=1}^k \Big\| E_j \Big( \sum_{j=1}^n u_i\Big)\Big\|
\le \delta \sqrt n\ .$$
\end{lem}

\begin{proof}
For $j\le k$ let $n_j = |\{ i:\supp u_i \cap E_j \ne \emptyset\}|$.
Thus $\sum_{j=1}^k n_j \le n$.
By Fact~2 and Cauchy-Schwarz,
\begin{equation*}
\frac1{2\sqrt k} \sum_{j=1}^k \Big\| E_j\Big( \sum_{i=1}^n u_i\Big)\Big\|
\le \frac{\sqrt 3}{2\sqrt k} \sum_{j=1}^k \sqrt{n_j}
\le \frac{\sqrt 3}{2\sqrt k} \Big( \sum_{j=1}^k n_j\Big)^{1/2}
\sqrt k \le \delta \sqrt n\ .
\end{equation*}
\end{proof}

\begin{proof}[Proof of Theorem~\ref{thm:3.2}]
Let $m\in\nat$.
We shall construct disjointly supported normalized vectors $(x_i)_1^m$ in
$Y$ so that $\|\sum_{i=1}^m x_i\| \le \sum_{i=0}^\infty \delta^i +1$.
By the unconditionality of $(e_i)$ this completes the proof.
Each $(x_i)$ will be a normalized average of certain basis vectors, the
number of which rapidly increases with $i$ and the supports
being uniformly mixed
(just as the 1 inch marks on a yardstick are uniformly separated by the
$\frac1{32}$ inch marks and so on).
To do this we need some notation.
We shall choose below rapidly increasing integers $q_1 < \cdots < q_m$.
Given these we define $p_i = \prod_{j=1}^i q_j$ for $i\le m$ and we then
choose natural numbers $r_{t_1,\ldots,t_i}$ for each $i\le m$ and
$t_j\le q_j$ for $j\le i$ so that $r_{t_1,\ldots,t_i} < r_{s_1,\ldots,s_j}$
whenever $(t_1,\ldots,t_i)$ is less than $(s_1,\ldots,s_j)$ lexicographically.
For example
$$r_1 < r_{2,1,4} < r_{2,2} < r_{2,2,3} < r_3\ .$$

We shall say that $(x_i)_{i=1}^m$ {\em corresponds\/} to $(q_1,\ldots,q_m)$
if for $i\le m$,
$$x_i = \frac2{\sqrt{p_i}} \sum_{j_1=1}^{q_1} \sum_{j_2=1}^{q_2}
\cdots \sum_{j_i=1}^{q_i} e_{r_{j_1,\ldots, j_i}}\ .$$
Since $(e_i)$ is 1-subsymmetric the particular choice of the
$r_{j_1,\ldots,j_i}$'s does not matter but their order does.
By Lemma~\ref{lem:3.4}, $\|x_i\| =1$ for $i\le m$.

Let $(\ep_n)_1^\infty$ be a sequence of positive numbers with
$\sum_{n=1}^\infty \ep_n <1$.
We shall prove by induction on $m$ that for every integer $q\ge 4$ there
exist integers $q<q_2 <\cdots < q_m$ so that for all integers $q_1$ with
$4\le q_1\le q$, if $(x_i)_1^m$ corresponds to $(q_1,\ldots,q_m)$ then
\begin{equation}\label{eq:3.2}
\Big\| \sum_{i=1}^m x_i\Big\| \le \sum_{i=0}^{m-1} \delta^i
+ \sum_{i=1}^m \ep_i\equiv M(m)\ .
\end{equation}

This is obvious for $m=1$ so assume it holds for some $m$.
Let $q\ge 4$.
Choose $d\in \nat$ so that
\begin{equation}\label{eq:3.3}
\sqrt{q} < \ep_{m+1} \sqrt{d}\ .
\end{equation}

Choose $n\in\nat$ so that
\begin{equation}\label{eq:3.4}
\frac{2dM(m)}{\sqrt{n}} < \ep_{m+1}\ .
\end{equation}

Let $q_2 = dn$.
By the inductive hypothesis for $q_0 \equiv qq_2$ we can find integers
$q_0 < q_3 <\cdots < q_{m+1}$ so that if $4\le s\le q_0$ and if $(y_i)_1^m$
is a sequence corresponding to $(s,q_3,q_4,\ldots,q_{m+1})$, then
$\|\sum_{i=1}^m y_i\| \le M(m)$.

Let $4\le q_1\le q$ and let $(x_i)_{i=1}^{m+1}$ correspond to
$(q_1,\ldots, q_{m+1})$.
There exists $k\ge2$ and $E_1 <\cdots < E_k$ so that
\begin{equation*}
\begin{split}
\Big\| \sum_{i=1}^{m+1} x_i\Big\|
& = \frac1{2\sqrt k} \sum_{j=1}^k \Big\| E_j\Big( \sum_{i=1}^{m+1} x_i\Big)
\Big\|\\
&\le \frac1{2\sqrt k} \sum_{j=1}^k \| E_j (x_1)\|
+ \frac1{2\sqrt k} \sum_{j=1}^k \Big\| E_j \Big( \sum_{i=2}^{m+1} x_i\Big)
\Big\| \ .
\end{split}
\end{equation*}
\medskip

\noindent {\bf Case 1.}
$k\ge d$.
$$\|x_1\|_1 = 2\sqrt{q_1} \le 2\sqrt{q} < 2 \ep_{m+1} \sqrt{d}$$
by \eqref{eq:3.3}.
Thus
$$\frac1{2\sqrt k} \sum_{j=1}^k \|E_j (x_1)\|
\le \frac1{2\sqrt k} \|x_1\|_1
< \frac{\sqrt d}{\sqrt k} \ep_{m+1} \le \ep_{m+1}\ .$$
Also
$$\frac1{2\sqrt k} \sum_{j=1}^k \Big\| E_j\Big( \sum_{i=2}^{m+1} x_i\Big)\Big\|
\le \Big\| \sum_{i=2}^{m+1} x_i\Big\|\ .$$
Now $(x_i)_{i=2}^{m+1}$ corresponds to the $m$-tuple
$(q_1q_2,q_3,\ldots,q_{m+1})$ and since $q_1q_2 \le qq_2 = q_0$ by the
inductive hypothesis,
$$\Big\| \sum_{i=2}^{m+1} x_i \Big\| \le M(m)\ .$$
Thus
$$\Big\|\sum_{i=1}^{m+1} x_i \Big\| \le M(m) + \ep_{m+1} < M(m+1)\ .$$
\medskip

\noindent {\bf Case 2.}
$k< d$.

For the $x_1$ term we use the estimate
$\frac1{2\sqrt k} \sum_{j=1}^k \|E_j (x_1)\| \le 1$.
To estimate the $\sum_{i=2}^{m+1} x_i$ term we write for
$2\le i\le m+1$, $x_i = \sum_{h=1}^n x_{i,h}$
where $(x_{i,h})_{h=1}^n$ is an identically distributed block basis.
Thus $|\supp (x_{i,h})| = p_i/n = q_1 dq_3 \ldots q_i$.

By Lemma~\ref{lem:3.4},
$\|x_{i,h}\| = \frac12 \sqrt{\frac{p_i}n} \frac2{\sqrt{p_i}}
= \frac1{\sqrt n}$ for $2\le h\le m+1$.
Thus for $h\le n$, $(\sqrt{n}\, x_{i,h})_{i=2}^{m+1}$ corresponds to
$(\frac{q_1q_2}{n},q_3,\ldots,q_{m+1})$ and since
$\frac{q_1q_2}{n} = q_1 d \le qq_2 = q_0$, by the inductive hypothesis we have
\begin{equation}\label{eq:3.5}
\Big\| \sum_{i=2}^{m+1} \sqrt{n}\, x_{i,h}\Big\| \le M(m)\ .
\end{equation}
Set $z_h = \sum_{i=2}^{m+1}  x_{i,h}$ for $h\le n$.
Then $(z_h)_{h=1}^n$ is an identically distributed block basis
of $(e_i)$ and hence
$$\|z_1\| = \cdots = \|z_n\| \equiv a \le \frac{M(m)}{\sqrt n}$$
by \eqref{eq:3.5}.

Since $E_1 <\cdots < E_k$, for $j\le k$ there are at most two $h$'s for
which $E_j \cap (\supp z_h) \ne\emptyset$ and
$E_{j'} \cap (\supp z_h) \ne \emptyset$ for some $j'\ne j$.

For $j\le k$ let
\begin{equation*}
\begin{split}
\tilde E_j = \cup & \{ \supp z_h : h\le n,\ \supp z_h\cap E_j\ne\emptyset\\
&\text{and }\supp z_h \cap E_{j'} = \emptyset\text{ if } j\ne j'\}\ .
\end{split}
\end{equation*}
We let $n_j$ be the cardinality of the set of such $h$'s.
Set $z= \sum_{i=2}^{m+1} x_i$.

Then since $a\le M(m)/\sqrt n$, $\|E_j z - \tilde E_jz\| \le
2M(m)/\sqrt n$.
Hence
$$\sum_{j=1}^k \|E_j z\| \le \frac{2M(m)k}{\sqrt n} + \sum_{j=1}^k
\|\tilde E_jz\|\ .$$
Now
$$\frac{2M(m)k}{\sqrt n} \le \frac{2M(m)d}{\sqrt n} <\ep_{m+1}
\ \text{ by \eqref{eq:3.4}}\ .$$
Also
$$\sum_{j=1}^k \|\tilde E_j z\|
\le 2a\sqrt{k}\, \delta \Big( \sum_{j=1}^k n_j\Big)^{1/2}
\le 2a \sqrt{k}\,\delta\sqrt{n}\ ,\ \text{ by Lemma~\ref{lem:3.5}.}$$
Thus $\frac1{2\sqrt k} \sum_{j=1}^k \|E_j (z)\| \le \ep_{m+1}
+ a\delta \sqrt{n}$.
Now $a\sqrt{n} \le M(m)$ and $\delta <1$ so this in turn is
$$< \delta \sum_{i=0}^{m-1} \delta^i + \sum_{i=1}^m \ep_i + \ep_{m+1}
= \sum_{i=1}^m \delta + \sum_{i=1}^{m+1} \ep_i\ .$$
{From} the $x_1$ estimate of $1= \delta^0$ we obtain
$\|\sum_{i=1}^{m+1} x_i\| \le M(m+1)$ in case~2.
\end{proof}

\begin{remarks}\label{rem:3.6}
A natural  question is whether $Y$ contains an a.s. subspace.
This is true by our next theorem.
The argument is motivated by arguments given in \cite{KOS}.
B.~Sari has also used a variation of these arguments to prove the
following results.
First we give some terminology from \cite{MMT}.

Let $(x_i)_{i=1}^\infty$ be a basis for a space $X$.
A normalized basic sequence $(d_i)_{i=1}^n \in \{X\}_n$ if $\forall\ \ep>0$
$\forall\ m_1\ \exists\ y_1\in \text{span}(x_i)_{i\ge m_1}$
$\forall\ m_2\ \exists\ y_2\in \text{span}(x_i)_{i\ge m_2}\cdots
\forall\ m_n \ \exists\ y_n\in \text{span}(x_i)_{i\ge m_n}$
so that $(y_i)_1^n$ is $(1+\ep)$-equivalent to $(d_i)_1^n$.
\end{remarks}

$X$ is {\em Asymptotic\/} $\ell_p$ if for some $K<\infty$
$\forall n \,\,
 \forall (d_i)_1^n\in \{X\}_n$, $(d_i)_{i=1}^n$ is $K$-equivalent
to the unit vector basis of $\ell_p^n$.
If $X$ is Asymptotic $\ell_p$ then $X$ contains an asymptotic $\ell_p$
basis, as defined above \cite{MMT}.

Note that we use a capital letter $A$ in the above definition in contrast
to the different notion asymptotic $\ell_p$ from the introduction. If $X$
is Asymptotic $\ell_p$, one can pass to a block basis which is
asymptotic $\ell_p$.

$X$ is {\em Asymptotically unconditional\/} if for some $K<\infty$
$\forall n \,\, \forall (d_i)_1^n\in \{X\}_n$, $(d_i)_{i=1}^n$ is
$K$-unconditional.
$A\buildrel K\over\sim B$ means $K^{-1} A\le B\le KA$.

\begin{THM}\label{thm:Sari}
\cite{Sa}
{\rm 1)}
Let $1<p<\infty$.
If $X$ is Asymptotically unconditional and for some $K<\infty$, for all
$m\le n\in \nat$ and for all disjointly supported normalized vectors
$(y_i)_{i=1}^m$ in $\text{span}(d_i)_{i=1}^n$,
$\|\sum_{i=1}^m y_i\| \buildrel K\over\sim m^{1/p}$,
then $X$ is Asymptotic $\ell_p$.

{\rm 2)} If $X$ is Asymptotically unconditional and for some $K<\infty$, for
all $n\in\nat$ and for all $(d_i)_{i=1}^n \in \{X\}_n$,
$\|\sum_{i=1}^n d_i\| \ge n/K$, then $X$ is Asymptotic $\ell_1$.
\end{THM}

\begin{thm}\label{thm:3.7}
Let $Z$ be a Banach space with a basis $(z_i)$.
Let $1\le p<\infty$ and $K<\infty$.
Assume that for all $(d_i)_{i=1}^n \in \{Z\}_n$,
$\|\sum_{i=1}^n d_i\| \buildrel K\over\sim n^{1/p}$.
Then every infinite dimensional subspace of $Z$ contains an asymptotic
$\ell_p$ basic sequence.
\end{thm}

\begin{cor}\label{cor:3.8}
Every infinite dimensional subspace of $Y= T_i (2;1/2)$ contains an
asymptotic $\ell_2$, hence a.s., subspace.
\end{cor}

Recall that a Banach space $X$ is minimal if every subspace of $X$
contains a further subspace isomorphic to $X$. Schlumprecht's space $S$
is minimal \cite{S2}.

\begin{cor}\label{cor:3.9}
The space $Y$ is not minimal.
\end{cor}

%Indeed, consider any subspace $Z$ of $Y$ which is generated by an
%asymptotic $\ell_2$ basic sequence $(z_i)$. If we assume that $Z$
%contains an isomorphic copy of $Y$, we can find in $Z$ a sequence
%$(u_i)$ equivalent to the standard basis $(e_i)$ of $Y$. Since $(e_i)$
%is subsymmetric, so is $(u_i)$. A subsequence of $(u_i)$ and thus the
%sequence itself is equivalent to a block basis of $(z_i)$, which is
%also asymptotic $\ell_2$. Therefore, $(e_i)$ is asymptotic $\ell_2$
%and so $Y$ is a.s., which is a contradiction. Alternatively, $(u_i)$ is
%a subsymmetric asymptotic $\ell_2$ sequence and thus it has to be
%equivalent to the u.v.b. of  $\ell_2$ and this is impossible \cite{CS}.
%Actually, what we use above is that $Z$ contains no subsymmetric basic
%sequences.

Indeed, Corollary~\ref{cor:3.8} yields that every subspace of $Y$
contains an asymptotic $\ell_2$ subspace $Z$. Since $Y$ does not
contain an isomorph of $\ell_2$ \cite{CS}, it follows that $Z$ cannot
contain a subsymmetric basic sequence and hence $Y$ does not embed
into $Z$.

We do not know if $Y$ contains a minimal subspace.

\begin{proof}[Proof of Theorem~\ref{thm:3.7}]
By standard perturbation arguments we need only show that every normalized
block basis $(x_i)$ of $(z_i)$ admits a further block basis which is
asymptotically $\ell_p$.
We may assume that $(z_i)$ is bimonotone, by renorming, and that $K>2$.
Furthermore, by passing to a block basis of $(x_i)$ we may assume
(\cite{MMT}, \cite{KOS}) that given $\ep_n\downarrow 0$ for $X= [(x_i)]$
\begin{equation}\label{eq:3.6}
\{X\}_n = \{W\}_n\text{ for all } n\in \nat \text{ and all block bases }
(w_i)\text{ of } (x_i);\ W= [(w_i)]\ .
\end{equation}
\begin{align}
&\text{for $n\in\nat$, if $(y_i)_{i=1}^n$ is a normalized block basis of
$(x_i)_{i=n}^\infty$,}
\label{eq:3.7}\\
&\text{then $(y_i)_1^n$ is $(1+\ep_n)$-equivalent to some $(d_i)_1^n\in
\{X\}_n$}.\notag
\end{align}
Thus, by increasing $K$, we may assume
\begin{align}
\text{For $n\in\nat$, if $(y_i)_{i=1}^n$ is }&\text{a normalized block basis
of $(x_i)_{i=n}^\infty$, then}
\label{eq:3.8}\\
&\Big\| \sum_{i=1}^n y_i\Big\| \buildrel K\over\sim n^{1/p}\ .
\notag
\end{align}

Fix $m\in\nat$ and $(d_i)_{i=1}^m \in \{X\}_m$.
We will prove that if $\sum_{i=1}^m |a_i|^p =1$ then
\begin{equation}\label{eq:3.9}
\frac12 \frac1{(2K)^2} \le \Big\| \sum_{i=1}^m a_i d_i\Big\| \le 2(2K)^2
\end{equation}
which will complete the proof of the theorem in view of our above remarks.

We choose $\ep,\delta,\delta' >0$ and $N' \in\nat$ to satisfy
\begin{gather}
0<\ep < \frac1{8m} \frac1{(2K)^{2p}}\label{eq:3.10}\\
\delta = \frac{\ep}{4Km}\label{eq:3.11}\\
N' > 2m\left( \frac{2K}{\delta}\right)^p\label{eq:3.12}\\
0<\delta' < \frac{\delta}{2KN'}\label{eq:3.13}
\end{gather}

Let $(w_i)_{i=1}^\infty$ be a normalized block basis of $(x_i)_{i=N'}^\infty$
so that for all $i$, if $w_i =\sum b_{i,j}x_j$ then $\sup_j |b_{i,j}|<\delta'$.
Such a $(w_i)$ exists by virtue of \eqref{eq:3.8}; one can let $w_i$ be
a suitably long average of $x_j$'s.
Given $\eta >0$ we can thus find, by \eqref{eq:3.6}, a normalized block
basis $(y_i)_{i=1}^m$ of $(w_i)$ which is $1+\eta$-equivalent to
$(d_i)_{i=1}^m$.
We will prove that \eqref{eq:3.9} holds for $(d_i)_1^m$ replaced by
$(y_i)_1^m$ and thus obtain \eqref{eq:3.9}.

Let $(a_i)_{i=1}^m \subseteq \real$ with $\sum_{i=1}^m |a_i|^p=1$.
{From} our construction we can write for $i\le m$, $a_iy_i =\sum_{j=1}^{n_i+1}
y_{i,j}$ where $n_i\ge0$, and $(y_{i,j})_{j=1}^{n_i+1}$ is a block
basis of $(x_i)_{i=N'}^\infty$, $\delta \le \|y_{i,j}\| <\delta +\delta'$
if $j\le n_i$ and $\|y_{i,n_i+1}\| <\delta$.
Set $\displaystyle N = \sum_{\substack{i\le m\\ |a_i|\ge \ep}} n_i$.
It follows from \eqref{eq:3.10} and \eqref{eq:3.11} that $N \ge 1$.
\begin{equation}\label{eq:3.14}
\text{If $i\le m$ and $|a_i|\ge \ep$ then } 1 \le n_i \le \left(\frac{2K}{\delta}
\right)^p\ .
\end{equation}
Indeed, suppose that $n_i \ge n_0 \equiv \llbrack (\frac{2K}{\delta})^p
\rrbrack +1$.
Then $1\ge \|a_i y_i\| \ge \frac{\delta}{K} n_0^{1/p} - n_0\delta'$ from
\eqref{eq:3.8}.
(We shrink each of $n_0$ successive $y_{ij}$'s to have norm exactly
$\delta$ at a loss of at most $\delta'$.
Note $n_0<N'$ so \eqref{eq:3.8} applies.) By \eqref{eq:3.10} we get
$n_i \ge 1$.

Now $n_0\delta' <\frac{\delta}{2K} n_0^{1/p}$ since this is equivalent
to $\delta' < \frac{\delta}{2Kn_0^{1/q}}$ (where $\frac1p + \frac1q=1$)
and we have
$$\delta' \buildrel {\eqref{eq:3.13}}\over <
\frac{\delta}{2KN'} \buildrel {\eqref{eq:3.12}}\over <
\frac{\delta \delta^p}{2K2m(2K)^p}
< \frac{\delta}{2Kn_0^{1/q}}$$
where the last inequality holds since $n_0^{1/q}\le n_0
< \frac{2m(2K)^p}{\delta^p}$.

Thus $1\ge \frac{\delta}{2K} n_0^{1/p}$ and so
$n_0 \le [\frac{2K}{\delta}]^p$, a contradiction.

Thus, by \eqref{eq:3.14} and \eqref{eq:3.12},
\begin{gather}
N\le m\left( \frac{2K}{\delta}\right)^p < N' \label{eq:3.15}\\
\text{If $i\le m$ and $|a_i| \ge \ep$ then }
\quad |a_i| = \| a_i y_i\| > \frac{\delta}{2K} n_i^{1/p} \label{eq:3.16}
\end{gather}

Indeed,  we can argue as in \eqref{eq:3.14} since $n_i \le N< N'$ to get
$$\|a_i y_i\| > \frac{\delta}K n_i^{1/p} - n_i \delta'\ .$$
Now $n_i\delta' < \frac{\delta}{2K} n_i^{1/p}$ is equivalent to
$\delta' < \frac{\delta}{2Kn_i^{1/q}}$.
But by \eqref{eq:3.13},
$\delta' < \frac{\delta}{2KN'} < \frac{\delta}{2Kn_i}$ and so
putting this together we obtain \eqref{eq:3.16}.
\begin{equation}\label{eq:3.17}
\text{If $i\le m$ and $|a_i|\ge \ep$ then }|a_i| = \|a_i y_i\|
< 2\delta Kn_i^{1/p}
\end{equation}

Again we have that $\|a_iy_i\| < \delta Kn_i^{1/p} + n_i\delta' +\delta$
by shrinking each $y_{i,j}$, $j\le n_i$, to have norm exactly $\delta$
at a cost of $\delta'$, and using \eqref{eq:3.8} and adding $\delta$ for
the term $\|y_{i,n_i+1}\|$.

We claim that $n_i\delta'+\delta < \delta Kn_i^{1/p}$.
Since $n_i\delta' +\delta <2\delta$, using $n_i < N'$ and \eqref{eq:3.13}
and since $2\delta <\delta Kn_i^{1/p}$ this yields \eqref{eq:3.17}.

Let $\displaystyle \mathop{{\sum}'}_{i=1}^m a_i y_i$ be the sum
$\displaystyle \sum_{\substack{|a_i|\ge \ep\\ 1\le i\le m}} a_i y_i$
and $\displaystyle \mathop{{\sum}''}_{i=1}^m a_i y_i$ be
$\displaystyle \sum_{\substack{|a_i|\ge \ep\\ 1\le i\le m}}
\sum_{j=1}^{n_i} y_{i,j}$.
We claim
\begin{equation}\label{eq:3.18}
\Big\| \mathop{{\sum}''}_{i=1}^m a_i y_i\Big\| \ge \frac{\delta}{2K}
N^{1/p}
\end{equation}

Indeed, by our now familiar method,
$$\Big\| \mathop{{\sum}''}_{i=1}^m a_i y_i\Big\| \ge \frac{\delta}{K}
N^{1/p} - N\delta'$$
but $N\delta' < \frac{\delta}{2K} N^{1/p}$ since
$$N\delta' < N'\delta' \buildrel {\eqref{eq:3.13}} \over <
\frac{\delta}{2K} < \frac{\delta}{2K} N^{1/p}\ .$$
Thus \eqref{eq:3.18} is proved.

{From} \eqref{eq:3.18} we have
\begin{equation*}
\begin{split}
\Big\| \mathop{{\sum}''}_{i=1}^m a_i y_i\Big\|^p >
\left( \frac{\delta}{2K}\right)^p N
&= \left(\frac{\delta}{2K}\right)^p
\mathop{{\sum}'}_{i=1}^m n_i
\buildrel {\eqref{eq:3.17}}\over >
\frac1{(2K)^p}
\mathop{{\sum}'}_{i=1}^m
\frac{|a_i|^p}{(2K)^p} \\
& = \frac1{(2K)^{2p}}
\mathop{{\sum}'}_{i=1}^m
|a_i|^p > \frac{(1-m\ep)}{(2K)^{2p}}
\end{split}
\end{equation*}
Thus
$$\Big\|\sum_{i=1}^m a_iy_i\Big\| >
\Big\| \mathop{{\sum}''}_{i=1}^m  a_iy_i\Big\| - m\ep -m\delta
> (1-m\ep)^{1/p}  \frac1{(2K)^2} -2m\ep
> \frac12 \frac1{(2K)^2} \ \text{ by \eqref{eq:3.10}}\ .$$
Next we show
\begin{equation}\label{eq:3.19}
\Big\| \mathop{{\sum}'}_{i=1}^m  a_iy_i\Big\| < 2\delta K N^{1/p}\ .
\end{equation}

As usual,
$\Big\| \mathop{{\sum}'}_{i=1}^m  a_iy_i\Big\| < \delta KN^{1/p}
+ N\delta' +m\delta$ and $N\delta' + m\delta < (m+1)\delta$.
We claim that $2m<K N^{1/p}$ which will complete the proof  of
\eqref{eq:3.19}.
First note that if $|a_i|\ge \ep$ then by \eqref{eq:3.17},
$n_i^{1/p} > \frac{\ep}{2\delta K} = 2m$ by \eqref{eq:3.11} which
yields \eqref{eq:3.19}.

Now
$$\Big\| \sum_{|a_i|<\ep} a_i y_i\| < m\ep < \frac12
\Big\| \mathop{{\sum}'}_{i=1}^m  a_iy_i\Big\| $$
since $\frac12
\Big\| {\displaystyle\mathop{{\sum}'}_{i=1}^m}  a_iy_i\Big\|
> \frac14 \frac1{(2K)^2} > m\ep$.

Thus by \eqref{eq:3.19},  $\|\sum_{i=1}^m a_i y_i\| < 3\delta KN^{1/p}$ so
\begin{equation*}
\begin{split}
&\Big\| \sum_{i=1}^m a_i y_i\Big\|^p < 3^p \delta^p K^p
\mathop{{\sum}'}_{i=1}^m n_i\\
&\qquad < 3^p K^p (2K)^p \mathop{{\sum}'}_{i=1}^m |a_i|^p\
\text{ by \eqref{eq:3.16}}\\
&\qquad < 3^p K^p (2K)^p
\end{split}
\end{equation*}
and this completes the proof.
\end{proof}

We do not  know if the modified space $V_{1/2,2}$ or if the modified
version of Schlumprecht's space $S$ \cite{S} are a.s.
Since their natural bases are symmetric, our arguments fail.
However we do have the next theorem.

\begin{thm}\label{thm:3.10}
Schlumprecht's space $S$ is not w.n.a.s.
\end{thm}

Recall that $S$ is the completion of $c_{00}$ under the norm which satisfies
the implicit equation
$$\|x\| = \max \Big\{ \|x\|_\infty, \sup \frac1{f(k)} \sum_{j=1}^k
\|E_j x\|\Big\}$$
where the ``sup'' is taken over all $k\ge 2$ and sets of integers
$E_1<\cdots < E_k$ with $f(k) = \log_2 (k+1)$.
It was shown in \cite{KL} that $c_0$ is finitely represented in $S$
and indeed our proof of Theorem~\ref{thm:3.2} was modeled after that
construction.
In \cite{M} an alternative ``partially nested'' construction was used.
We shall follow the proof of Theorem~3.1 in \cite{M} and use his
notation as we sketch the proof of Theorem~\ref{thm:3.10}.

Take $p=1$, $q=\infty$, $\theta_k = 1/f(k)$, $n_k=K$,
$1/p_k=1-\log_k (1/\theta_k) = 1-\frac{\ln f(k)}{\ln k}$ and
$1/q_k = \frac{\ln f(k)}{\ln k}$.

In the course of the proof, if $n$ is arbitrary, $k_0$ is chosen so that
$1/f(k_0)\le 1/n$ and then $k_1$ is chosen with $f(k_0)/f(k_1) \le 1/n$.
If we take $k_0$ so that $f(k_0)$ has ``order $n$,'' we need $k_1$ to
satisfy (basically)
\begin{equation}\label{eq:3.20}
f(k_1) \ge n^2\ .
\end{equation}
Also the proof in \cite{M} requires $k_0^{1/q_{k_1}} \le 2$ which
essentially transforms into
\begin{equation}\label{eq:3.21}
\frac{ n\ln f(k_1)}{\ln k_1} \le 1\ .
\end{equation}

\begin{proof}[Proof of Theorem~\ref{thm:3.10} (sketch)]
For an arbitrary integer $n$ let $d= k_1$ satisfy \eqref{eq:3.20} and
\eqref{eq:3.21}.
To prove $S$ is not w.n.a.s.\ we need only check the condition for
$m\equiv nd$.

In order to apply Theorem~\ref{thm:3.1} \cite{M} we choose a rapidly
increasing sequence of integers $m\ll q_1 \ll \cdots \ll q_n$.
We shall employ $n$ different normalized distributions of elements of $S$.
Precisely, for $j\le n$ set
$$v_j = \frac{f(q_j)}{q_j} \sum_{i=1}^{q_j}  e_i$$
where $(e_i)$ is the unit vector basis of $S$.
We then successively repeat each $v_j$ $d$-times so that altogether we
have $m=nd$ vectors to serve as distributions.
We shall compare the norms of two different permutations of $m$ of a
block sequence of $m$ vectors with these distributions.
Since $(e_i)$ is 1-subsymmetric this will show that $S$ is not w.n.a.s.

Let $(y_{i,j})_{i=1,\, j=1}^{d,n}$ be a block basis of $(e_i)$ in
lexicographic order with $y_{i,j}$ equal to $v_j$ in distribution.
For $j\le n$ let $u_j = \sum_{i=1}^d y_{i,j}$.
Then
$$\|u_j\| = \frac{f(q_j)dq_j}{q_jf(dq_j)}\approx d\ .$$
The proof of Theorem~\ref{thm:3.1} \cite{M} yields a constant $C$,
independent of $m$, such that $(d^{-1}u_j)_{j=1}^n$ is $C$-equivalent
to the unit vector basis of $\ell_\infty^n$.
Thus
\begin{equation}\label{eq:3.22}
\Big\|\sum_{i=1}^d \sum_{j=1}^n  y_{i,j}\Big\|
= \Big\| \sum_{j=1}^n u_j\Big\| \le Cd = \frac{Cm}{n}\ .
\end{equation}

We next consider a different order.
Let $(z_{i,j})_{i=1,\, j=1}^{d,n}$ be a block basis of $(e_i)$,
ordered with respect to the lexicographic order $(j,i)$, the reverse
coordinates, with $z_{i,j}$ equal in distribution to $v_j$.
Set $w_j = \sum_{i=1}^d z_{i,j}$ so that $(w_j)_{j=1}^d$ is a block
basis of $(e_i)$.
As before
\begin{equation}\label{eq:3.23}
\Big\| \sum_{i=1}^d \sum_{j=1}^n z_{i,j}\Big\|
= \Big\| \sum_{j=1}^n w_j\Big\|
\ge \frac{n}{f(n)} \frac{d}2 = \frac12 \frac{m}{f(n)}\ .
\end{equation}

Since $\frac{m}{f(n)}/\frac{m}n = \frac{n}{f(n)}$ we have from
\eqref{eq:3.22} and \eqref{eq:3.23} that $S$ is not w.n.a.s.
We obtain more precisely that the constant $C_m$ for $m$ sequences is
at least of the order $\frac{\sqrt{\ln (m)}}{\ln(\ln (m))}$.
\renewcommand{\qed}{}
\end{proof}

\begin{remarks}\label{rem:3.11}
Clearly the constant $C_m$, the w.n.a.s. constant for $m$ normalized
weakly null sequences in $S$, satisfies $C_m\le f(m)$.
We  can show that $C_2 = f(2)$ with a different construction.
Pei-Kee Lin \cite{L} has pointed out that one can adjust the proof in
\cite{KL} to obtain a sequence in $S$ whose spreading model is isometric
to $\ell_1$.
Indeed in \cite{KL} it was shown that for every $\ep>0$ there exists a
rapidly increasing sequence of integers, $(p_k)_{k=1}^\infty$, so that
if $u_j = \frac{f(p_j)}{p_j} \sum_{i=1}^{p_j} e_i$ then for all $n$
there exist disjointly supported vectors $(v_j)_{j=1}^n$ in $S$, each of
the same distribution as $u_j$, with $\|\sum_{j=1}^n v_j\|\le 1+\ep$.
One can choose scalars $(a_k)_{k=1}^\infty \subseteq (0,1)$ converging
to 1 so that $\|z_n\|<1$ and $\|z_n\|\to 1$ if $z_n = \sum_{j=1}^n a_jv_j$.
It follows that any block basis $(x_n^1)$ with distribution $x_n^1=$
distribution $z_n$ has spreading model 1-equivalent to the unit vector
basis of $\ell_1$.
If we let $(x_n^2)$ be a block basis with distribution $x_n^2=$
distribution $u_n$, then the w.n.a.s.\ constant for these two
sequences is $f(2)$.
Finally we note that since $S$ is minimal \cite{S2}, no subspace of $S$ is
w.n.a.s.
\end{remarks}

%%%%%%%%%%%%%%  44444
\section{Tsirelson's space is not iteration stable}

No good criterion is known which forces a Banach space $X$ to be isomorphic
to a stable space.
Attempting to find such a criterion, H.~Rosenthal has asked if it might be
true that every asymptotically symmetric and iteration stable space $X$
is isomorphic to a stable space.
In this regard he asked us if $T$ is iteration stable.
We show that it is not in this section.
First we give the relevant definitions.

\begin{defin}\label{def:4.1}
A sequence $(x_n)$ in a Banach space $X$ is {\em type determining\/}
if for all $x\in X$,
$$\lim_{n\to\infty} \|x+x_n\| \text{ exists.}$$
\end{defin}

\begin{defin}[H. Rosenthal]\label{def:4.2}
A Banach space $X$ is {\em iteration stable\/} if for all type determining
sequences $(x_n)$ and $(y_n)$ in $X$,
$$\lim_{n\to\infty} \lim_{m\to\infty} \|x_n + y_m\|\ \text{ exists.}$$
\end{defin}

Iteration stability is another softening of the definition of stability.

\begin{thm}\label{thm:4.3}
Tsirelson's space $T$ is not iteration stable.
\end{thm}

We recall that $T$ is the completion of $c_{00}$ under the norm satisfying
the implicit equation
$$\|x\| = \max \bigg( \|x\|_\infty, \sup \frac12 \sum_{i=1}^n \|Ex\|\bigg)$$
where the ``sup'' is taken over all $n\ge2$ and $n\le E_1<\cdots < E_n$.
The unit vector basis $(e_i)$ is a normalized 1-unconditional basis for the
reflexive space $T$.
For $k\ge 2$ we set
$$\|x\|_k = \sup \frac12 \sum_{i=1}^k \|E_ix\|$$
where the ``sup'' is taken over all $k\le E_1 <\cdots < E_k$.

\begin{proof}[Proof of Theorem~\ref{thm:4.3}]
For $n\ge 2$ let $z_n = \frac2{n^2} \sum_{i=1}^{n^2} e_{n^3+i}$ so that
$\|z_n\|=1$ and $z_n$ is an $\ell_1^n$-average with constant~1.
Precisely, $z_n = \frac1n \sum_{j=1}^n z_{n,j}$ where $z_{n,j} =
\frac2n \sum_{i=1}^n e_{n^3 + (j-1) n+i}$ for $j\le n$.
$(z_{n,j})_{j=1}^n$ is 1-equivalent to the unit vector basis of $\ell_1^n$.
A standard calculation yields
\begin{equation}\label{eq:4.1}
\lim_n \|z_n\|_k = \frac12\ \text{ for all }\ k\ge 2\ .
\end{equation}

For $n\in\nat$ we set
$$x_n = \begin{cases} e_{n^3} + \frac14 e_{n^3+1}\ ,&\text{$n$ odd}\\
\noalign{\vskip4pt}
e_{n^3}  +\frac14 z_n\ ,&\text{$n$ even.}
\end{cases}$$
We let $y_m = \frac12 \sum_{i=1}^4 e_{m+i}$ for $m\in\nat$.

$(y_n)$ is normalized and clearly type determining.
Also $\|x_n\| = \|x_n\|_\infty=1$ for all $n$.
Moreover for $k\ge 2$, $n\in\nat$,
\begin{equation}\label{eq:4.2}
\|x_n\|_k \le \frac12 \Big(1+\frac12\Big) = \frac34\ .
\end{equation}

Note that
$$\lim_{n\text{ odd}} \lim_m \|x_n + y_m\| = \frac12 \Big( 1+\frac14 +2\Big)
= \frac{13}8\ .$$
Also
$$\lim_{n\text{ even}} \lim_m \|x_n + y_m\| = \frac12 \Big( 1+\frac12 +2\Big)
= \frac{14}8\ .$$

We shall show that $(x_n)$ is type determining and this will complete the
proof.
To do this we verify by induction on $\ell$ that $\lim_{n\to\infty} \|x+x_n\|$
exists for all $x\in\text{span}(e_i)$ with $|\supp x| \le \ell$.
The case $\ell = 0$ holds since $\|x_n\|=1$ for all $n$ so assume it holds
for $\ell$ and let $|\supp x| = \ell +1$.

We may assume that $\supp x<\supp x_n$ and we define
$$S(x+x_n) = \sup_{k\le \max (\supp x)} \sup \frac12
\sum_{i=1}^k \|E_i (x+x_n)\|$$
where the second ``sup'' is taken over all $2\le k\le E_1 <\cdots < E_k$ with
$E_1 \cap \supp x\ne \emptyset$.
We will show that
\begin{equation}\label{eq:4.3}
\lim_{n\to\infty} S(x+x_n)\ \text{ exists.}
\end{equation}

Observe that for any $n$,
$$\|x+x_n\| = \max \{\|x+x_n\|_\infty, S(x+x_n)\}\ .$$
Indeed if $k\le E_1 <\cdots < E_k$ with $\supp x<E_1$, then
by \eqref{eq:4.2}
$$\frac12 \sum_{i=1}^k \|E_i (x+x_n)\|
\le \|x_n\|_k \le \frac34 < \|x_n\|_\infty\ .$$
Thus \eqref{eq:4.3} will complete the proof.

To prove \eqref{eq:4.3} we will show that we can balance the various types
of splitting of $x+x_n$ by $(E_i)_1^k$ between $n$ odd and $n$ even.
There is no need to consider splittings such that the last set $E_k$
intersects $\supp x$ since that case is covered by the induction hypothesis.
Also any values obtained for $S(x+x_n)$ resulting from splittings for $n$
odd can be treated by a similar splitting for $n$ even, treating
$\frac14 z_n$ as a single element.
We need to show that $\varlimsup_{n\text{ even}} S(x+x_n)$ cannot be
bigger than $\lim_{n\text{ odd}} S(x+x_n)$.

Consider a splitting $(E_i)_1^k$ giving rise to $S(x+x_n)$ for $n$ even.
\medskip

\noindent {\bf Case 1.}
There exists $i_0 <k$ such that $\max E_{i_0} = n^3$.
Thus the support of $z_n$ is split into $k-i_0$ intervals, $k\le\max \supp x$.
By \eqref{eq:4.1} this presents no problem as we can obtain  same values
for $n$ odd (as $n\to\infty$).
\medskip

\noindent {\bf Case 2.}
There exists $i_0 <k$ such that $n^3 \in \supp E_{i_0}$ and
$E_{i_0} \cap \supp z_n \ne \emptyset$.
By the triangle inequality,
\begin{equation*}
\begin{split}
&\frac12 \sum_{i=1}^k \|E_i (x+x_n)\|
\le \frac12 \sum_{i=1}^{i_0-1} \|E_i x\| \\
&\qquad +\frac12 \| E'_{i_0} (x+x_n)\|
+ \frac12 \|E''_{i_0} x_n\|
+ \frac12 \sum_{i=i_0+1}^k \|E_i x_n\|
\end{split}
\end{equation*}
where $E_{i_0} = E'_{i_0} \cup E''_{i_0}$, $E'_{i_0} < E''_{i_0}$ and
$\max E'_{i_0} = n^3$.
The new splitting might not be admissible but $j= k-i_0+1\le \max\supp x$
and again by \eqref{eq:4.1}, the value for the new splitting behaves in the
limit as an admissible splitting $(E_1,\ldots, E_{i_0-1},E'_{i_0},E)$
where $E = E''_{i_0}\cup \bigcup_{i=i_0+1}^k E_i$.
This last splitting can be mimicked to yield the same value for $n$ odd.
\end{proof}

\begin{remark}\label{rem:4.4}
Our proof depended upon using vectors whose norm is given by
$\|\cdot\|_\infty$.
In regard to Rosenthal's question it is of interest to determine if $T$
contains iteration stable subspaces.
Accordingly we have the following partial result.
\end{remark}

\begin{prop}\label{prop:4.5}
Let $X\subseteq T$ be an infinite dimensional subspace so that for all
$0\ne x\in X$, $\|x\|\ne \|x\|_\infty$.
Then if $(x_n),(y_n)\subseteq X$ are weakly null sequences,
$\lim_n \|x_n\|$ exists and $(y_n)$ is type determining, we have
\begin{equation}\label{eq:4.4}
\lim_{n\to\infty} \lim_{m\to\infty} \|x_n+ y_m\| \ \text{ exists.}
\end{equation}
\end{prop}

\begin{proof}
We may assume that $\|x_n\|=1$ for all $n$.
Also, to prove \eqref{eq:4.4} we can freely pass to subsequences of $(y_n)$
and so we may assume that for all $k$
\begin{equation}\label{eq:4.5}
\lim_m \|y_m\|_k = \lambda (k)\ \text{ exists.}
\end{equation}
By perturbing we may assume that $(x_n)$ and $(y_m)$ are block bases
of $(e_i)$.

Since $\|y_m\|_k \le \|y_m\|_{k+1}$ when $\supp y_m>k$, it follows that
$\lambda (k) \le \lambda (k+1) \le \lim_m \|y_m\|$ for all $k$.
Let $\lambda (k) \uparrow \lambda$.
We prove
\begin{equation}\label{eq:4.6}
\lim_n \lim_m \|x_n + y_m\| = (1+\lambda) \vee \lim_m \|y_m\|
\end{equation}
which yields \eqref{eq:4.4}.

Let $\ep >0$ and choose $k\in\nat$ with $\lambda (k)> \lambda-\ep$.
Choose $\bar t \in \nat$ so that if $t\ge \bar t$ and
$\sum_{i=1}^t |a_i|=1$ then there exists $F\subseteq \{1,\ldots,t\}$ with
$|F|=k$ and $\sum_{i\in F} |a_i| <\ep$.
Let $n$ be such that $\supp x_n >\bar t$.
Choose $m_0$ so that for $m\ge m_0$, $\big|\| y_m\|_k - \lambda (k)\big|<\ep$.
We will prove first that
\begin{equation}\label{eq:4.7}
\varliminf_m \|x_n + y_m\| > (1+\lambda -2\ep) \vee \lim_m \|y_m\|\ .
\end{equation}

Let $1= \|x_n\| = \frac12 \sum_{i=1}^t \|E_i x_n\|$ for some $t\ge \bar t$
and $t\le E_1 <\cdots < E_t$.
Let $m\ge m_0$ so that $\supp x_n < \supp y_m$ and choose
$\min (\supp y_m) \le F_1 <\cdots < F_k$ with $\|y_m\|_k = \frac12
\sum_{i=1}^k \|F_i y_m\|$.
By deleting the smallest $k$ terms from $\frac12 \sum_{i=1}^t \|E_i x_n\|$
and replacing them by $\frac12 \sum_{i=1}^k \| F_i y_m\|$ we obtain
$$\|x_n + y_m\| \ge 1-\ep + \|y_m\|_k > 1+\lambda -2\ep$$
and \eqref{eq:4.7} follows.

We next prove that
\begin{equation}\label{eq:4.8}
\varlimsup_m \|x_n + y_m\| < (1+\lambda +\ep) \vee \lim_m \|y_m\|
\end{equation}
which will complete the proof of \eqref{eq:4.6}.

Let $\supp x_n < \supp y_m$ and let $\|x_n + y_m\| = \frac12 \sum_{i=1}^t
\|G_i(x_n + y_m)\|$ where $t\ge \bar t$ and $t\le G_1 <\cdots < G_t$.
Suppose there does not exist $i_0$ with $G_{i_0} x_n\ne0$ and
$G_{i_0} y_m \ne 0$.
Then
$$\|x_n + y_m\| \le \|x_n\| + \|y_m\|_{k(m)}$$
where $k(m) \le t$.
If this occurred for infinitely many $m$'s we would obtain
$\lim_m \|x_n + y_m\| \le (1+\lambda) \vee \lim_m \|y_m\|$.
Suppose such a $G_{i_0}$ exists.
Split $G_{i_0} = G'_{i_0} \cup G''_{i_0}$ with $G'_{i_0} < G''_{i_0}$
and $G_{i_0} x_n = G'_{i_0} x_n$.
We now have $t+1$ sets in the sum
$$\frac12 \sum_{i=1}^{i_0-1} \|G_i x_n\| + \|G'_{i_0} x_n\|
+ \|G''_{i_0} y_m\| + \sum_{i=i_0+1}^t \|G_i y_m\|\ .$$
If we throw out the smallest term we have $t$ sets after coordinate $t$
and so for some $\ell (m)\le \max \supp x_n$,
$$\|x_n + y_m \| \le \|x_n\| + \|y_m\|_{\ell (m)} + \ep\ .$$
Again, letting $m\to\infty$, we obtain \eqref{eq:4.8}.
\end{proof}

\end{document}